\documentclass[12pt]{amsart}
\usepackage{amssymb,latexsym}
\newdimen\AAdi%
\newbox\AAbo%
\def\AAk#1#2{\s_etbox\AAbo=\hbox{#2}\AAdi=\wd\AAbo\kern#1\AAdi{}}%
\def\AAr#1#2#3{\s_etbox\AAbo=\hbox{#2}\AAdi=\ht\AAbo\raise#1\AAdi\hbox{#3}}%

\font\tenmsb=msbm10 at 12pt
\font\sevenmsb=msbm7 at 8pt
\font\fivemsb=msbm5 at 6pt
\newfam\msbfam
\textfont\msbfam=\tenmsb
\scriptfont\msbfam=\sevenmsb
\scriptscriptfont\msbfam=\fivemsb
\def\Bbb#1{{\tenmsb\fam\msbfam#1}}
\def\R{\Bbb R}
\def\C{\Bbb C}
\def\B{\Bbb B}
\def\N{\Bbb N}
\def\Q{\Bbb Q}
\def\Z{\Bbb Z}
\def\PP{\Bbb P}
\def\EE{\Bbb E}
\def\F{\Bbb F}
\def\G{\Bbb G}
\def\H{\Bbb H}
\def\S{\Bbb S}
\begin{document}
\newtheorem{thm}{Theorem}
\newtheorem{lem}{Lemma}
\newtheorem{cor}{Corollary}
\newtheorem{rem}{Remark}
\newtheorem{pro}{Proposition}
\newtheorem{defi}{Definition}
\newcommand{\noi}{\noindent}
\newcommand{\dis}{\displaystyle}
\newcommand{\mint}{-\!\!\!\!\!\!\int}
\newcommand{\ba}{\begin{array}}
\newcommand{\ea}{\end{array}}
\def \bx{\hspace{2.5mm}\rule{2.5mm}{2.5mm}} \def \vs{\vspace*{0.2cm}} 
\def\hs{\hspace*{0.6cm}}
\def \ds{\displaystyle}
\def \p{\partial}
\def \O{\Omega}
\def \o{\omega}
\def \b{\beta}
\def \m{\mu}
\def \l{\lambda}
\def \ul{u_\lambda}
\def \D{\Delta}
\def \d{\delta}
\def \e{\varepsilon}
\def \a{\alpha}
\def \tf{\widetilde{f}}
\def\cqfd{%
\mbox{ }%
\nolinebreak%
\hfill%
\rule{2mm} {2mm}%
\medbreak%
\par%
}
\def \pr {\noindent {\it Proof.} }
\def \rmk {\noindent {\it Remark} }
\def \esp {\hspace{4mm}}
\def \dsp {\hspace{2mm}}
\def \s_sp {\hspace{1mm}}
\def \u{u_+^{p^*}}
\def \ui{(u_+)^{p^*+1}}
\def \ul{(u^k)_+^{p^*}}
\def \energy{\int_{\R^n}\u }
\def \s{\sigma}
\def \sk{\s_k}
\def \mo{\mu_k}
\def \I{{\cal I}}
\def \J{{\cal J}}
\def \K{{\cal K}}
\def \OM{\overline{M}}
\def\fk{{{\cal F}}_k}
\def\M1{{{\cal M}}_1}
\def\n{\nabla}
\def\uuu{{\n ^2 u+du\otimes du-\frac {|\n u|^2} 2 g_0}}
\def\sku{\s_k^{1/k}\left(\uuu\right)}
\def\L{\Lambda}
\def\vvv{{\frac{\n ^2 v} v -\frac {|\n v|^2} {2v^2} g_0+S_{g_0}}}
\def\vvs{{\frac{\n ^2 \widetilde v} {\widetilde v}
 -\frac {|\n \widetilde v|^2} {2\widetilde v^2} g_{S^n}+S_{g_{S^n}}}}
\def\skv{\sk\left(\vvv\right)}
\def\tr{\hbox{tr}}
\def\pO{\partial \Omega}
\def\dist{\hbox{dist}}
\def\RR{\Bbb R}\def\R{\Bbb R}
\def\C{\Bbb C}
\def\B{\Bbb B}
\def\N{\Bbb N}
\def\Q{\Bbb Q}
\def\Z{\Bbb Z}
\def\PP{\Bbb P}
\def\EE{\Bbb E}
\def\F{{\mathcal F}}
\def\G{\Bbb G}
\def\H{\Bbb H}
\def\qed{\cqfd}
\def\NMI{{Newton-MacLaurin inequality}}
\title[Schouten tensor]{Some properties of the Schouten tensor and applications to conformal geometry}
\author{Pengfei Guan}
\address{Department of Mathematics\\
 McMaster University\\
Hamilton, Ont. L8S 4K1, Canada.\\
Fax: (905)522-0935 }
\email{guan@math.mcmaster.ca}
\thanks{Research of the first author was supported in part by 
NSERC Grant OGP-0046732.}
\author{Jeff Viaclovsky}
\address{Department of Mathematics, MIT, Cambridge, MA, USA}
\email{jeffv@math.mit.edu}
\author{Guofang Wang}
\address{Max-Planck-Institute for Mathematics in
the Sciences\\ Inselstr. 22-26, 04103 Leipzig, Germany}
\email{gwang@mis.mpg.de}
\begin{abstract}
The Riemannian curvature tensor decomposes into a conformally 
invariant part, the Weyl tensor, and a non-conformally 
invariant part, the Schouten tensor. A study of the 
$k$th elementary symmetric function of the eigenvalues 
of the Schouten tensor was initiated in \cite{Jeff1},
and a natural condition to impose is that the 
eigenvalues of the Schouten tensor are  
in a certain cone,  $\Gamma_k^+$. 
We prove that this eigenvalue condition for $k \geq n/2$ implies 
that the Ricci curvature is positive.
We then consider some applications to the 
locally conformally flat case, in particular, 
to extremal metrics of $\sigma_k$-curvature functionals 
and conformal quermassintegral inequalities, 
using the results from \cite{GW}.
\end{abstract}
\maketitle
\section{Introduction}
Let $(M^n,g)$ be an $n$-dimensional  
Riemannian manifold, $n \geq 3$, and let the Ricci 
tensor and scalar curvature be denoted by  
$Ric$ and $R$, respectively. 
We define the Schouten tensor
\begin{align*}
A_g= \frac{1}{n-2}\left( Ric - \frac{1}{2(n-1)}Rg \right).
\end{align*}
There is a decomposition formula (see \cite{Besse}):
\begin{equation}\label{decomp}
\mbox{Riem} = A_g \odot g + {\mathcal W}_g,
\end{equation}
where ${\mathcal W}_g$ is the Weyl tensor of $g$, and $\odot$ denotes the Kulkarni-Nomizu product (see \cite{Besse}). 
As Weyl tensor is conformally
invariant, to study the deformation of conformal metric, we only
need to understand the Schouten tensor.  A study of $k$-th elementary 
symmetric functions of the Schouten tensor was initiated in \cite{Jeff1},  
it was reduced to certain fully nonlinear
Yamabe type equations. In order to apply elliptic theory of fully 
nonlinear equations, one often restricts Schouten tensor to be
in certain cone $\Gamma^+_k$ defined as follows (according to 
G{\.a}rding \cite{Garding}).

\begin{defi}
Let $(\lambda_1, \cdots, \lambda_n) \in \mathbf{R}^n$. Let 
$\sigma_k$ denote the $k$th elementary symmetric function
$$\sigma_k(\lambda_1, \cdots, \lambda_n) = \sum_{i_1 < \cdots < i_k}
\lambda_{i_1} \cdots \lambda_{i_k},$$ and we let $$\Gamma_k^+ =
\mbox{component of } \{\sigma_k > 0\} \mbox{ containing} \quad (1,\cdots,1).$$ 
Let $\bar \Gamma_k^+$ denote the closure of $\Gamma_k^+$. If 
$(M,g)$ is a Riemannian manifold, and $x \in M$, we say $g$ has 
positive (nonnegative resp.) $\Gamma_k$-curvature at $x$ if its Schouten
tensor $A_g \in \Gamma^+_k$ ($\bar \Gamma^+_k$ resp.) at $x$. In this
case, we also say
$g \in \Gamma^+_k$ ($\bar \Gamma^+_k$ resp.) at $x$. 
\end{defi}

We note that positive $\Gamma_1$-curvature is equivalent to positive scalar curvature, 
and the condition of positive $\Gamma_k$-curvature has some
geometric and topological consequences for the manifold $M$. For example,
when $(M,g)$ is locally conformally flat with positive
$\Gamma_1$-curvature, then $\pi_i(M)=0, \forall 1<i\leq \frac n2$ by a 
result of Schoen-Yau \cite{SY}. 
In this note, we will prove that positive $\Gamma_k$-curvature
for any $k\geq \frac n2$ implies positive Ricci
curvature.

\begin{thm} \label{kbig} Let $(M,g)$ be a Riemannian manifold
and $x \in M$, if $g$ has positive (nonnegative resp.) $\Gamma_k$-curvature 
at $x$ for some $k\ge
 n/2$. Then its Ricci curvature is positive (nonnegative resp.) at $x$.
Moreover,  if  $\Gamma_k$-curvature is nonnegative for some  $k>1$, 
then
\[ Ric_g \ge \frac{2k-n}{2n(k-1)} R_g \cdot g .\]
In particular if $k\ge \frac n2$, 
\[ Ric_g \ge \frac{(2k-n)(n-1)}{(k-1)} 
{\binom{n}{k}}^{-\frac 1k} 
\sigma^{\frac 1k}_k(A_g) \cdot g .\]
\end{thm}
 
\noindent{\it Remark.} Theorem \ref{kbig} is not true for $k=1$.
Namely the condition of positive scalar curvature
 gives no restriction on the lower bound of Ricci curvature .
 
\begin{cor} \label{cor1} Let $(M^n,g)$ be a compact,
locally conformally flat manifold with nonnegative
$\Gamma_k$-curvature  everywhere for some $k\ge n/2$. Then $(M,g)$ is 
conformally equivalent to
either a  space form or 
a finite quotient of a Riemannian $\S^{n-1}(c)\times \S^1$
for some constant $c>0$ and $k=n/2$. Especially, if
$g\in \Gamma_k^+$, then $(M,g)$ is conformally equivalent to
a spherical space form.
\end{cor}
  
 When $n=3,4$, the result in Theorem~\ref{kbig} was already observed in \cite{GV} and \cite{CGY}. Theorem~\ref{kbig} and Corollary~\ref{cor1}
will be proved in the next section. 

We will also consider the equation 
\begin{align}
\label{eqn1}
\sigma_k (A_{\tilde{g}}) = constant,
\end{align}
for conformal metrics $\tilde{g} = e^{-2u}g$.
This equation was studied in \cite{Jeff1}, where 
it was shown that when $k \neq n/2$, (\ref{eqn1}) is the 
conformal Euler-Lagrange equation of the functional 
\begin{equation}\label{functional}
{ \mathcal F}_k(g)=vol(g)^{-\frac{n-2k} n}\int_M \sk(g)\, dvol(g),
\end{equation}
when $k=1, 2$ or for $k>2$ when $M$ is locally conformally 
flat. 
We remark that in the even dimensional locally 
conformally flat case, $\mathcal{F}_{n/2}$ is a conformal 
invariant.  Moreover, it is a multiple of the Euler 
characteristic, see \cite{Jeff1}.  

 This problem was further studied studied in \cite{GW},
where the following conformal flow was considered:
\begin{align*}
\frac{d}{dt} g &= - ( \log \sigma_k(g) - \log r_k(g)) \cdot g,\\
g(0) &= g_0,
\end{align*}
where 
\[\log r_k = \frac{1}{Vol(g)} \int_M \log 
\sigma_k(g) dvol(g).\]
Global existence with uniform $C^{1,1}$ a priori bounds
of the flow was proved in \cite{GW}. It was also proved that
 for $k \neq n/2$ the flow is sequentially convergent in $C^{1, \alpha}$ 
to a $C^{\infty}$ solution of $\sigma_k = constant$. 
Also, if $k < n/2$, then $\mathcal{F}_k$ is decreasing along 
the flow, and if $k > n/2$, then $\mathcal{F}_k$ is 
increasing along the flow. We remark that the existence 
result for equation (\ref{eqn1}) has been obtained 
independently in \cite{Yanyan} in the locally conformally 
flat case for all $k$. 

 In Section 3, we will consider global 
properties of the functional $\mathcal{F}_k$, and 
give conditions for the existence of a global 
extremizer. We will also derive some 
conformal quermassintegral inequalities, 
which are analogous to the classical quermassintegral 
inequalities in convex geometry. 

\section{Curvature restriction}
We first state a proposition which describes some important properties of
the sets $\Gamma_k^{+}$.
\begin{pro}
\label{coneprop}
(i) Each set
$\Gamma_k^{+}$ is an open convex cone with vertex at the
origin, and we have the following sequence of inclusions
$$\Gamma_n^{+} \subset \Gamma_{n-1}^{+} \subset \cdots
\subset \Gamma_{1}^{+}.$$
(ii) For any $\L=(\l_1,\cdots, \l_n) \in \Gamma^+_k$ ($\bar \Gamma^+_k$ resp.), 
$\forall 1\le i \le n$, let $$(\L|i)=(\l_1,\cdots,\l_{i-1},\l_{i+1},\cdots,\l_n),$$
then $(\L|i) \in \Gamma^+_{k-1}$ ($\bar \Gamma^+_{k-1}$ resp.). In particular,
$$\Gamma_{n-1}^{+} \subset V_{n-1}^{+} = \{ (\lambda_1, \cdots, \lambda_n) \in \mathbf{R}^n
: \lambda_i + \lambda_j > 0, i \neq j \}.$$
\end{pro}
The proof of this proposition is standard following from \cite{Garding}.

\medskip

Our main results are the consequences of the following two lemmas.
In this note, we assume that $k>1$.
 
 \begin{lem}\label{lem1} Let 
 $\L=(\l_1,\l_2,\cdots, \l_{n-1}, \l_n)\in \R^n$, and define
 \[A_\L=\L-\frac{\sum_{i=1}^n\l_i}{2(n-1)}(1,1,\cdots, 1).\]
If $A_\L \in \bar \Gamma_k^+$, then
\begin{equation}\label{eq0}
\min_{i=1, \cdots, n} \l_i\ge \frac{(2k-n)}{2n(k-1)} \sum_{i=1}^n\l_i.\end{equation}
In particular if $k\ge \frac n2$, 
\[ \min_{i=1, \cdots, n} \l_i \ge \frac{(2k-n)(n-1)}{(n-2)(k-1)}
{ \binom{n}{k}}^{-\frac 1k} 
\sigma^{\frac 1k}_k(A_\L).\]
\end{lem}

\pr  We first note that,
for any non-zero vector $A=(a_1,\cdots,a_n) \in \bar \Gamma^+_2$
implies $\sigma_1(A)>0$. This can be proved as follow. As $A \in \bar \Gamma^+_2$, 
$\sigma_1(A)\ge 0$. If $\sigma_1(A)=0$,  
there must be $a_i>0$ for some $i$ since $A$ is a non-zero vector. 
We may assume $a_n>0$.
Let $(A|n)=(a_1,\cdots,a_{n-1})$, we have $\sigma_1(A|n)\ge 0$
by Proposition~\ref{coneprop}. This would give
$\sigma_1(A)=\sigma_1(A|n)+a_n>0$, a contradiction.

Now without loss of generality, we may
assume that $\L$ is not a zero vector. By the assumption
$A_\L \in \bar \Gamma_k^+$ for $k\ge 2$, so we have $\sum_{i=1}^{n} \l_i>0$. 

 Define 
\[\L_0=(1,1, \cdots, 1,\d_k)\in \R^{n-1}\times\R\]
and we have $A_{\L_0}=(a,\cdots,a, b)$, where
\[\d_k=\frac {(2k-n)(n-1)}{2nk-2k-n},\]
\[a=1-\frac{n-1+\d_k}{2(n-1)}, \quad b=\d_k-\frac{n-1+\d_k}{2(n-1)}\] 
so that
 \begin{equation}\label{eq1}
 \s_k(A_{\L_0})=0 \quad \text{ and }\s_{j}(A_{\L_0})>0 \text{ for } j\le k-1.\end{equation}
 It is clear that $\d_k<1$ and so that $a>b$.
Since (\ref{eq0}) is invariant under the transformation
$\L$ to $s\L$ for $s>0$, we may assume that $\sum_{i=1}^n\l_i= \tr(\L_0)=n-1+\d_k$
and $\l_n=\min_{i=1, \cdots, n} \l_i$. We write 
$$A_\L=(a_1,\cdots,a_n).$$
We claim that
\begin{equation}\label{eq2}
\l_n\ge \d_k.\end{equation}
This is equivalent to show
\begin{equation}\label{eq300}
a_n\ge b.
\end{equation}
Assume by contradiction that $a_n<b$.
We consider $\L_t=t\L_0+(1-t)\L$ and
\[A_t:=A_{\L_t}=tA_{\L_0}+(1-t)A_{\L}=
((1-t)a+ta_1,\cdots, (1-t)a+ta_{n-1}, (1-t)b+ta_n).\]
By the convexity of the cone $\Gamma_k^+$ (see Proposition 1),
we know 
\[ A_t \in \bar \Gamma_k^+, \quad \text{ for any }t\in (0,1].\]
Especially, $f(t):=\s_k(A_t)\ge 0$ for any $t\in [0,1].$ 
By the definition of $\d_k$, $f(0)=0$.

For any $i$ and any vector $V=(v_1,\cdots,v_n)$,
we denote $(V|i)=(v_1,\cdots,v_{i-1},v_{i+1},\cdots,v_n)$  
be the vector with the $i$-th
component removed.
Now we compute the derivative of $f$ at $0$
\[
f'(0)= \sum_{i=1}^{n-1} (a_i-a)\s_{k-1}(A_0|i)+ (a_n-b)\s_{k-1}(A_0|n).\]
Since $(A_0|i)=(A_0|1)$ for $i\le n-1$ and $\sum _{i=1}^na_i =(n-1)a+b$, 
we have
\[
f'(0)  =\ds\vs (a_n-b)
(\s_{k-1}(A_0|n)-\s_{k-1}(A_0|1))
<0,\]
for $\s_{k-1}(A_0|n)-\s_{k-1}(A_0|1) >0$. (Recall that $b<a$.)
This is a contradiction, hence $\l_n\ge \d_k$. It follows that
\[\min_{i=1, \cdots, n}\l_i \ge \d_k =\frac{2k-n}{2n(k-1)}\sum_{i=1}^n\l_i.\]
Finally, the last inequality in the lemma follows from the Newton-MacLaurin
inequality. \qed

\noindent
{\it{Remark.}} It is clear from the above proof that 
the constant in Lemma \ref{lem1} is optimal. \\

We next consider the case $A_\L\in \bar \Gamma_{\frac n2}^+$.
\begin{lem}\label{lem2} Let $k=n/2$ and $\L=(\l_1,\cdots,\l_n)\in \R^n$ with $A_\L\in
\bar \Gamma_k^+$. 
Then either $\l_i>0$ for any $i$ or 
\[\L=(\l,\l,\cdots,\l,0)
\]
up to a permutation. If the second case is true, 
then we must have $\sigma_{\frac n2}(A_\L)=0$.
\end{lem}
\pr By Lemma \ref{lem1}, to prove the Lemma we
only need to check that for $\L=(\l_1, \cdots,\l_{n-1},0)$ with
 $A_\L \in \bar \Gamma_k^+$, 
 \[\l_i=\l_j, \quad \forall i,j=1,2,\cdots, 2k-1.\]
We use  the same idea as in the proof of the previous Lemma. 
Without loss of generality, we may
assume that $\L$ is not a zero vector. By the assumption
$A_\L \in \bar \Gamma_k^+$ for $k\ge 2$,  we have $\sum_{i=1}^{n-1} \l_i>0$. 
Hence we may assume that $\sum_{i=1}^{n-1}\l_i=n-1$.
Define 
\[\L_0=(1,1, \cdots, 1,0)\in \R^n\]
 It is easy to check that
\begin{equation}\label{1}
A_{\L_0} \in \Gamma_{k-1}^+ \quad \text{and} \quad
\s_k(A_{\L_0})=0.\end{equation}
That is, $ A_{\L_0} \in \bar \Gamma_{k}^+$.
If $\l$'s are not all the same, we have 
$$\sum_{i=1}^{n-1} (\l_i-1)=0,$$
and
$$\sum_{i=1}^{n-1} (\l_i-1)^2>0.$$
Now consider $\L_t=t\L_0+(1-t)\L$ and
\[A_t:=A_{\L_t}=tA_{\L_0}+(1-t)A_{\L}=
(\frac 12+t(\l_1-1),\cdots, \frac 12+t(\l_{n-1}-1),-\frac 12).\]
From the assumption that $A\in \bar \Gamma_k^+$, (\ref{1}) and the convexity
of $\bar \Gamma_k^+$, we have
\begin{equation}\label{2}
 A_t \in \bar \Gamma_k^+  \quad \text{for } t>0.\end{equation}

For any $i\neq j$ and any vector $A$,
we denote $(A|ij)$ be the vector with
the $i$-th and $j$-th components removed. Let 
$\widetilde \Lambda=(\frac 12,\cdots,
\frac 12, -\frac 12)$ be $n-1$-vector, $\Lambda^*=(\frac 12,\cdots,
\frac 12, -\frac 12)$ be $n-2$-vector. It is clear
that $\forall i \neq j, \quad i,j \le n-1$, 
$$\sigma_{k-1} (A_{0}|i)=\sigma_{k-1}(\widetilde \Lambda)>0 ,$$
$$\sigma_{k-2}(A_{0}|ij)=\sigma_{k-2}(\Lambda^*)>0.$$
 Now we expand $f(t)=\s_k(A_t)$ at $t=0$. By (\ref{1}),
 $f(0)=0$. We compute
 \[\begin{array}{rcl}
 f'(0)& =& \ds\vs \sum_{i=1}^{n-1}(\l_i-1)\s_{k-1}(A_0|i) \\
 &=&
 \ds \s_{k-1}(\widetilde \Lambda)\sum_{i=1}^{n-1}(\l_i-1)=0\end{array}\]
 and
  \[\begin{array}{rcl}
 f''(0)& =& \ds\vs \sum_{i\not = j}(\l_i-1)(\l_j-1)\sigma_{k-2}(A_{0}|ij) \\
&=&
 \ds\vs \s_{k-2}(\Lambda^*)\sum_{i\not = j}(\l_i-1)(\l_j-1)\\
&=& \ds -\s_{k-2}(\Lambda^*)\sum_{i=1}^{n-1}(\l_i-1)^2<0,\\
 \end{array}\]
 for $\sigma_{k-2}(A_{0}|ij)=\s_{k-2}(\Lambda^*)>0$ for any $i\not =j$ 
 and $\sum_{i \neq j} ( \lambda_j - 1)= (1 - \lambda_i)$.
 Hence $f(t)<0$ for small $t>0$, which contradicts (\ref{2}). 
 \qed

\noindent
{\it{Remark.}} From the proof of Lemma~\ref{lem2}, 
there is a constant $C>0$ depending
only on $n$ and $\frac{\sigma_{\frac n2}^{\frac 2n}(A_\L)}{\sigma_1(A_\L)}$
such that 
$$\min_i \l_i \ge C \sigma_{\frac n2}^{\frac 2n}(A_\L).$$ \\

\medskip

\noindent{\it Proof of Theorem \ref{kbig}.} Theorem~\ref{kbig} follows 
directly from Lemmas \ref{lem1} and \ref{lem2}.  \qed
\medskip
\begin{cor} \label{product} Let $(M,g)$ is a $n$-dimensional
Riemannian manifold and $k\ge n/2$,
and let $N=M\times \S^1$ be the product manifold.
Then $N$ does not have positive $\Gamma_k$-curvature. If $N$ has
nonnegative $\Gamma_k$-curvature, then 
$(M,g)$ is an Einstein manifold, and there are two cases:
either $k=n/2$ or $k>n/2$ and $(M,g)$ is a torus.
\end{cor}

\pr This follows from Lemmas \ref{lem1} and \ref{lem2}. \qed

\medskip

\noindent{\it Proof of Corollary~\ref{cor1}.} From Theorem~\ref{kbig},
we know that the Ricci curvature $Ric_g$ is nonnegative. Now
 we deform it by the Yamabe flow considered by Hamilton, Ye \cite{Ye} 
 and Chow \cite{Chow} to obtain a conformal
 metric $\widetilde g$ of constant scalar curvature. The
 Ricci curvature $Ric_{\widetilde g}$ is nonnegative, for the Yamabe flow
 preserves the non-negativity of Ricci curvature, see \cite{Chow}.
 Now by a classification result given in \cite{Tani, Chen}, we have
 $(M,\widetilde g)$ is isometric to either a  space form
 or a finite quotient of a Riemannian $\S^{n-1}(c)\times \S^1$
for some constant $c>0$. In the latter case, it is clear
 that $k=n/2$, otherwise it can not have nonnegative $\Gamma_k$-curvature.
  \qed

\bigskip

Next, we will prove that if $M$ is locally conformally flat
with positive $\Gamma_{n-1}$-curvature,
then $g$ has positive sectional curvature.
If $M$ is locally conformally flat, then by (\ref{decomp}) we may decompose 
the full curvature tensor as  
$$\mbox{Riem} = A_g \odot g.$$

\begin{pro}
\label{sectional}
Assume that $n=3$, or that $M$ is locally conformally flat. Then
Schouten tensor $A_g \in V_{n-1}^{+}$ if and only if $g$ has positive
sectional curvature. 
\end{pro}

\pr
Let  $\pi$ be any $2$-plane in $T_p(N)$, and let $X$,$Y$ be
an orthonormal basis of $\pi$. We have that
\begin{align*}
K(\sigma)&= \mbox{Riem}(X,Y,X,Y) = A_g \odot g (X,Y,X,Y)\\
&= A_g(X,X)g(Y,Y)- A_g(Y,X) g(X,Y) + A_g(Y,Y)g(X,X) - A_g(X,Y)g(Y,X)\\
&= A_g(X,X) + A_g(Y,Y).
\end{align*}
From this it follows that
$$ {\underset{\sigma \in T_pN}{\mbox{min}}}
K(\sigma) = \lambda_1 + \lambda_2,$$
where $\lambda_1$ and $\lambda_2$ are the smallest eigenvalues
of $A_g$ at $p$. 
\qed

\begin{cor}
If $(M,g)$ is locally conformally flat
with positive $\Gamma_{n-1}$-curvature,
then $g$ has positive sectional curvature. 
\end{cor}
\pr
This follows easily from Propositions \ref{coneprop}
and \ref{sectional}.
\qed

\section{Extremal metrics of $\sigma_k$-curvature functionals}
We next consider some properties of the functionals $\mathcal{F}_k$ 
associated to $\sigma_k$. These functionals were introduced 
and discussed in \cite{Jeff1}, see also \cite{GW}.
Further variational properties in connection to $3$-dimensional 
geometry were studied in \cite{GV}. 

\medskip

We recall that ${\mathcal F}_k$ is defined by
\[
{ \mathcal F}_k(g)=vol(g)^{-\frac{n-2k} n}\int_M \sk(g)\, dvol(g).
\]
We denote ${\mathcal C}_k=\{g \in [g_0] | g \in \Gamma_k^+\}$,
where $[g_0]$ is the conformal class of $g_0$.

We now apply our results to show that if $g_0 \in \Gamma^+_{\frac n2}$, then
there is an extremal metric $g_e$ which minimizes 
${\mathcal F}_{m}$ for $m < n/2$, and if
$m > n/2$, there is an extremal metric $g_e$ which 
maximizes ${\mathcal F}_{m}$.

\medskip

\begin{pro}\label{extremal}
Suppose $(M,g_0)$ is locally conformally flat and
$g_0 \in \Gamma^+_k$ for some $k \geq \frac n2$, then 
$\forall  m < \frac n2$ there is an extremal
metric $g^m_e \in [g_0]$ such that
\begin{equation}\label{extrem1}
\inf_{g \in {\mathcal C}_m}  { \mathcal F}_m(g)={ \mathcal F}_m(g^m_e), 
\end{equation}
and $\forall  m > \frac n2$, there is extremal metric $g^m_e \in [g_0]$
such that
\begin{equation}\label{extrem2}
\sup_{g \in {\mathcal C}_m}  { \mathcal F}_k(g)={ \mathcal F}_k(g^m_e), 
\end{equation}
In fact, any solution to $\sigma_m(g)=constant$ is an extremal metric.
\end{pro}
\pr 
First by Corollary \ref{cor1}, $(M, g_0)$ is conformal to 
a spherical  space form. For any $g\in {\mathcal C}_m$,
from \cite{GW} we know there is a conformal metric $\widetilde g$ in ${\mathcal C}_m$
satisfying that $\sigma_m(\widetilde g)$  is constant 
and
\begin{itemize}
\item[(a).] if $m>n/2$, then $\F_m(g) \le \F_m(\widetilde g)$.
\item[(b).] if $m<n/2$, then $\F_m(g) \ge \F_m(\widetilde g)$.
\end{itemize}
A classification result of \cite{Jeff1}, \cite{Jeff2} which 
is analogous to a result of Obata for the 
scalar curvature, shows that $\widetilde g$ has
constant sectional curvature.
Therefore $\widetilde{g}$ is the unique critical metric 
unless $M$ is conformally equivalent to $\S^n$, 
in which case any critical metric is 
the image of the standard metric 
under a conformal diffeomorphism. 
The clearly implies the conclusion of the 
Proposition.
\qed

  Next we consider the case $k < n/2$. We have 

\begin{pro}\label{extremal2}
Suppose $(M,g_0)$ is locally conformally flat and
$g_0 \in \Gamma^+_k$ for some $k < \frac n2$.
Suppose furthermore that for any fixed $C >0$, the
space of solutions to the equation $\sigma_k = C$ 
is compact, with a bound independent of the 
constant $C$. Then there is an extremal
metric $g^k_e \in [g_0]$ such that
\[\inf_{g \in {\mathcal C}_k}  { \mathcal F}_k(g)={ \mathcal F}_k(g^k_e).\]
\end{pro}
\pr
From the compactness assumption, 
there exists a critical metric $g^k_e$ which 
has least energy. If the functional assumed a 
value strictly lower than ${ \mathcal F}_k(g^k_e)$,
then by \cite{GW}, the flow would decrease to another 
solution of $\sigma_k = constant$, which is a contradiction 
since $g^k_e$ has minimal energy. 
\qed

\medskip
  We conclude with conformal quermassintegral 
inequalities, which were speculated in \cite{GW}, 
and verified there for some special 
cases when $(M,g)$ is locally conformally flat and 
$g \in \Gamma^+_{\frac n2-1}$ 
or $g \in \Gamma^+_{\frac n2+1}$ using the flow method.
In the case of $k=2, n=4$, the inequality was proved in \cite{Gursky}
without the locally conformally flat assumption. 

\begin{pro}\label{afi}
Suppose $(M,g_0)$ is locally conformally flat and
$g_0 \in \Gamma^+_k$ for some $k \geq \frac n2$, then for any $1\le
 l<\frac n2 \le k\le n$ there is a constant $C(k,l,n)>0$, such that
for any $g \in [g_0]$ and $g \in \Gamma^+_k$
\begin{equation}\label{cqm}
({\mathcal F}_k(g))^{1/k} \le C(k,l,n)({ \mathcal F}_l(g))^{1/l},
\end{equation}
with equality if and only if $(M,g)$ is a spherical space form.
\end{pro}

\pr By Proposition \ref{extremal}, we have a conformal metric $g_e$ 
of constant sectional curvature satisfies
such that
\[\inf_{g \in {\mathcal C}_l}  { \mathcal F}_l(g)={ \mathcal F}_l(g_e)\]
and
 \[\sup_{g \in {\mathcal C}_k}  { \mathcal F}_k(g)={ \mathcal F}_k(g_e).\]
Hence, we have for any $g\in \Gamma_k^+$
\[ \begin{array}{rcl}
\ds\vs \frac{(\F_k(g))^{1/k}} {(\F_l(g))^{1/l}}
&\le& \ds\frac{(\F_k(g_e))^{1/k}} {(\F_l(g_e))^{1/l}}\\
&=&\ds \frac{(l!(n-l)!)^{1/l}}{(k!(n-k)!)^{1/k}}.\end{array}\]
When  the equality holds, $g$ is an extremal of $\F_l$, hence
a metric of constant sectional curvature by \cite{Jeff1}.
\qed

\end{document}